\documentclass[12pt]{article}
\usepackage{amsfonts,amsmath,amssymb,amscd}
\usepackage{shadow}
\usepackage{graphicx}
\usepackage{graphicx}
\usepackage{color}

\newtheorem{theo}{Theorem}[section]

\newtheorem{coro}[theo]{Corollary}

\def\qed{\hfill \rule{4pt}{7pt}}

\def\proof{\noindent {\it{Proof.} \hskip 2pt}}

\begin{document}

\begin{center}
{\large\bf Labeled Ballot Paths and the Springer Numbers}
\end{center}

\begin{center}
William Y.C. Chen$^1$, Neil J.Y. Fan$^2$, Jeffrey Y.T. Jia$^3$

Center for Combinatorics, LPMC-TJKLC\\
Nankai University, Tianjin 300071, P.R. China

$^1$chen@nankai.edu.cn, $^2$fjy@cfc.nankai.edu.cn,
$^3$jyt@cfc.nankai.edu.cn.

\end{center}

\vskip 6mm \noindent {\bf Abstract.}
The Springer numbers are defined
in connection with the irreducible root systems of type $B_n$,
which also arise as
the generalized Euler and class numbers introduced by Shanks.
Combinatorial interpretations of the Springer numbers have been found  by
Purtill  in terms of   Andr\'e
signed permutations,
 and by Arnol'd   in terms of  snakes of type
$B_n$.
We introduce the inversion code of a snake of type $B_n$
and establish a bijection
 between labeled ballot paths of length $n$
 and snakes of type $B_n$.  Moreover, we obtain the bivariate
 generating function for the number $B(n,k)$
 of labeled ballot paths starting at $(0,0)$
 and ending at  $(n,k)$.
 Using our bijection,
 we  find a statistic $\alpha$
 such that the number of snakes $\pi$ of type $B_n$ with
  $\alpha(\pi)=k$ equals $B(n,k)$. We also show
   that our bijection specializes  to
   a bijection between
   labeled Dyck paths of length $2n$ and
  alternating permutations on $[2n]$.



\section{Introduction}

 The Springer numbers  are
 introduced  by Springer \cite{Springer}
 in the study of irreducible root system of type  $B_n$.
  Let $S_n$ denote the $n$-th Springer number. The sequence $\{S_n\}_{n\geq 0}$
  is listed  as  entry $A001586$ in OEIS  \cite{Sloane}. The first few values of $S_n$ are
  \[1,1,3,11,57,361,2763,24611, \ldots.\]
 To be more specific,  $S_n$ can be defined as follows.
 Let $V$ be a real vector space, $R$ be a  root system  of type $B_n$ in $V$,
 and $W$ be the Weyl group of $R$. It is known that for a fixed simple root set $S$ of
$R$, any $\alpha\in R$ is either a positive or a negative linear
combination of elements of $S$, denoted by
 $\alpha>0$ or $\alpha<0$. For a subset $I\subset S$, let
 $\sigma(I,S)$ denote the number of elements $w\in W$ such that $w\alpha>0$
  for any $\alpha\in I$ and $w\alpha<0$ for any $\alpha\in S-I$.
  Then the Springer number  $S_n$ can be  defined as
  the maximum value of $\sigma(I,S)$ among $I\subset S$.
  Springer derived the following generating
function,
\begin{align}\label{gf}
\sum\limits_{n\geq 0} S_n \frac{x^n}{n!}=\frac{1}{\cos x -\sin x}.
\end{align}

On the other hand, Hoffman \cite{Hoff} pointed out that the Springer
numbers also arise as the generalized Euler and class numbers
$s_{m,n}$ $(n\geq 0)$ for $m=2$, where the numbers $s_{m,n}$  are introduced by
Shanks \cite{Shanks} based on the Dirichlet series
 \begin{align*}
L_{m}(s)=\sum\limits_{k=0}^{\infty}
\left(\frac{-m}{2k+1}\right)\frac{1}{(2k+1)^{s}}.
\end{align*}
Note that the above notation $(-m/(2k+1))$ is the Jacobi symbol. To be precise,
the generalized Euler and class numbers $s_{2,n}$ are defined by
\[
 s_{2,n}=\left\{
       \begin{array}{ll}
          c_{2,\frac{n}{2}}, &\mbox{ if $n$ is even;}\\[5pt]
          d_{2,\frac{n+1}{2}}, &\mbox{ if $n$ is odd,}
       \end{array}
       \right.
\]
where the numbers $c_{2,n}$ and $d_{2,n}$ are
  given by
\begin{equation}
c_{2,n}=\frac{(2n)!}{\sqrt 2}\left(\frac{\pi}{4}\right)^{-2n-1}L_2(2n+1),\nonumber
\end{equation}
\begin{equation}
d_{2,n}=\frac{(2n-1)!}{\sqrt 2}\left(\frac{\pi}{4}\right)^{-2n}L_{-2}(2n).\nonumber
\end{equation}

According to the following recurrence relations for $c_{2,n}$ and $d_{2,n}$ derived by
 Shanks \cite{Shanks},
\[
\sum\limits_{i=0}^{n}(-4)^{i}{2n \choose 2i}c_{2,n-i}=(-1)^n,
\]
\[
\sum\limits_{i=0}^{n-1}(-4)^{i}{2n-1 \choose 2i}d_{2,n-i}=(-1)^{n-1},
\]
one sees that the numbers $s_{2,n}$ are integers.
In fact, the above recurrence relations lead to the following formulas
\begin{align*}
\sum\limits_{n\geq 0}c_{2,n}\frac{x^{2n}}{(2n)!}&=\sec  2x \cos x ,\\[5pt]
\sum\limits_{n\geq 1}d_{2,n}\frac{x^{2n-1}}{(2n-1)!}&=\sec  2x \sin
x .
\end{align*}

Shanks raised the question of finding  combinatorial interpretations
for the Euler and class numbers $s_{m,n}$. For $m=2$,    $s_{2,n}$
is the $n$-th Springer number. Purtill \cite{Purtill} gave  an
interpretation of the Springer numbers  in terms of the Andr\'e
signed permutations on $[n]=\{1,2,\ldots, n\}$.  Arnol'd \cite{Arno}
found another interpretation of the Springer numbers in terms of
snakes of type $B_n$. Recall that
 a snake of type $B_n$ is an alternating signed permutation
 $\pi=\pi_1\pi_2\ldots\pi_n$ on $[n]$ such that
 \begin{equation}\label{ud}
  0<\pi_1>\pi_2<\pi_3>\pi_4<\cdots \pi_n.
  \end{equation}
 For example, $1\bar 3 2$ is a snake of type $B_3$. Intuitively, a signed
permutation on $[n]$ can be viewed as  an ordinary permutation
$\pi_1\pi_2\cdots\pi_n$ with some elements associated with minus
signs. An element $i$ with a minus sign is often written as
$\bar{i}$. The above alternating or up-down condition (\ref{ud})
is based on the following natural order:
\[
\bar n<\ldots<\bar 1<1<\ldots<n.
\]
 Arnol'd \cite{Arno} proved that the Springer number $S_n$ equals the number of snakes of
type $B_n$. Hoffman \cite{Hoff} showed that the exponential
generating function for the number of snakes of type $B_n$ also
equals the right hand side of (\ref{gf}), that is, the generating
function of the Springer numbers. Recently, Chen, Fan and Jia
\cite{CFJ} obtained a formula for the generating function of
$s_{m,n}$ for arbitrary $m$, which in principle leads to a
combinatorial interpretation of the numbers $s_{m,n}$ in terms of
 alternating augmented  $m$-signed permutations. Note that for $m=2$,
alternating augmented 2-signed permutations are exactly  snakes of
type $B_n$.

The objective of this paper is to give a combinatorial
interpretation for the Springer numbers in terms of labeled ballot
paths. In fact, we shall introduce the inversion code of a snake
of type $B_n$. By using the inversion code, we construct a
bijection between the set of snakes of type $B_n$ and the set of
labeled ballot
 paths of length $n$.
 Let
$B(n,k)$ denote the number of labeled ballot paths starting at
$(0,0)$ and ending at $(n,k)$. Then the numbers $B(n,k)$ can be
viewed as a refinement of the Springer numbers. Using the
recurrence relation of $B(n,k)$,
we obtain the generating function
for $B(n,k)$ for any given $n$.

Using our bijection, we find a statistic $\alpha$  on
 snakes
of type $B_n$ such that the number of snakes $\pi$ of type $B_n$
with $\alpha(\pi)= k$ equals $B(n,k)$. A labeled ballot path that
eventually returns to the $x$-axis is called a labeled Dyck path.
When $k=0$, $B(2n,0)$ is the number of labeled Dyck paths of
length $2n$. We find that $B(2n,0)$ and the number $E_{2n}$ of
alternating permutations on $[2n]$ have the same generating
function, and we show that our bijection for labeled ballot paths
and snakes of type $B_n$ reduces to a bijection between labeled
Dyck paths and alternating permutations.

The paper is organized as follows. In Section 2, we give
descriptions of  the map $\Phi$ from snakes of type $B_n$ to
 labeled ballot paths of
length $n$, and the map $\Psi$ from labeled ballot paths of length
$n$ to snakes of type $B_n$. In Section 3, we prove that the maps
$\Phi$ and $\Psi$ are well defined, and they are inverses
 of each
other. The last section is devoted to the bivariate
generating
function for the numbers $B(n,k)$  and
 the classification of snakes of type $B_n$
 in accordance with the numbers $B(n,k)$. We also show that
 the map $\Psi$ restricted to labeled Dyck paths serves
 as a combinatorial interpretation of the fact that
$B(2n,0)$ equals $E_{2n}$.

\section{The bijection }

In this section, we define a class of labeled ballot paths and
establish a bijection between such labeled ballot paths of length
$n$ and snakes of type $B_n$. Recall that a  ballot path of length
$n$  is a lattice path with $n$ steps from the origin  consisting
of up steps $u=(1,1)$ and down steps $d=(1,-1)$ that do not go
below the $x$-axis. As a special case,
 a Dyck path  is a ballot
path of length $2n$ that ends at  the $x$-axis. A ballot path is
also called a partial Dyck path  \cite{Chen}. The height of a step
of a ballot path is defined to be
the smaller $y$-coordinate of
its endpoints. By a labeled ballot
 path we mean a ballot path for
which each step is endowed with a nonnegative integer that is less
than or equal to its height. A labeled
 ballot path $P=p_1p_2\cdots
p_n$ for which step $p_i$ is labeled by $w_i$ is denoted by
$(P;W)$, where $W=w_1w_2\cdots w_n$.

\begin{figure}[h,t]
\begin{center}
\begin{picture}(130,50)
\setlength{\unitlength}{1mm}

\put(0,0){\circle*{1.5}}\put(0,0){\line(1,1){5}}\put(1,1.8){\small0}
\put(5,5){\circle*{1.5}}\put(5,5){\line(1,1){5}}\put(6,6.8){\small1}
\put(10,10){\circle*{1.5}}\put(10,10){\line(1,1){5}}\put(11,11.8){\small1}
\put(15,15){\circle*{1.5}}\put(15,15){\line(1,-1){5}}\put(17,12){\small0}
\put(20,10){\circle*{1.5}}\put(20,10){\line(1,-1){5}}\put(22,7.4){\small1}
\put(25,5){\circle*{1.5}}\put(25,5){\line(1,1){5}}\put(26.5,7.5){\small1}
\put(30,10){\circle*{1.5}}\put(30,10){\line(1,1){5}}\put(31.5,12.5){\small2}
\put(35,15){\circle*{1.5}}

\end{picture}
\caption{A labeled ballot path $(uuudduu;0110112)$ of length
$7$.}\label{ff1}
\end{center}
\end{figure}
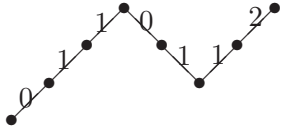

For example, for a ballot path $P=uuudduu$,  there are 216
labelings. Figure \ref{ff1} gives a labeling of the ballot path $P$.

For $n=3$, there are 3 ballot paths $P_1=uuu, P_2=uud$ and
$P_3=udu$. There are 6 labelings for $P_1$, 4 labelings
for $P_2$
and 1 labeling for $P_3$.
On the other hand, there are 11 snakes
of type $B_3$ as listed below:
\[1\bar23,1\bar32,1\bar3\bar2,213,2\bar13,2\bar31,2\bar3\bar1,312,3\bar12,3\bar21,3\bar2\bar1.\]

In order to establish a bijection between
 ballot paths of length
$n$ and snakes of type $B_n$, we introduce the
inversion code of a snake $\pi$ of type $B_n$.  Write $\pi=\pi_1\cdots\pi_{n}$. We define
$c_i(\pi)$
 as follows
\[
c_i(\pi)=\left\{
              \begin{array}{ll}
                \#\{(\pi_{2k},\pi_{2k+1})| 1\leq k\leq  (n-1)/2, i<2k,\pi_{2k}<\pi_i<\pi_{2k+1}\}, & \hbox{if $n$ is odd;} \\[5pt]
                \#\{(\pi_{2k-1},\pi_{2k})| 1\leq k\leq n/2, i<2k-1,\pi_{2k}<\pi_i<\pi_{2k-1}\}, & \hbox{if $n$ is even.}
              \end{array}
            \right.
\]
The
 sequence $(c_1(\pi), c_2(\pi), \ldots, c_{n}(\pi))$, denoted $c(\pi)$, is called the inversion code of $\pi$. For
example, let $n=7$ and $\pi=3\bar5214\bar76$. Then
 the inversion code  of $\pi$ is $(2,1,2,1,1,0,0)$.
For $n=8$ and  $\pi=538\bar2\bar1\bar476$, the inversion code of
$\pi$ is  $(1,1,0,1,0,0,0,0)$.

We are now ready to describe the map $\Phi$ from   a snake
$\pi=\pi_1\pi_2\cdots\pi_n$ of type $B_n$ to a labeled ballot path
$(P;W)=(p_1p_2\cdots p_n; w_1w_2\cdots w_n)$. Suppose that $p_1,
p_2, \ldots, p_{k-1}$ and their labels $w_1,w_2,\ldots,w_{k-1}$
have been determined, we proceed to demonstrate
 how to determine
$p_k$ and its label $w_k$. If we were in Step 1, namely, for
$k=1$, we would locate the element $n$ or $\bar{n}$ in $\pi$, and
would assume that $\pi_i=n$ or $\bar{n}$. Suppose that we are
 in Step $k$. Now we look for the
 element $n-k+1$ or $\overline{n-k+1}$ in $\pi$. Here are two
 cases.

\noindent Case 1.  Assume that $\pi_i=n-k+1$. If $i$ is odd, then
set $p_k=u$; if $i$ is even, then set $p_k=d$. Set $w_k=c_i(\pi)$.

\noindent Case 2. Assume that $\pi_i=\overline{n-k+1}$. If $i$ is
odd, then set $p_k=d$; if $i$ is even, then set $p_k=u$. Set
$w_k=h_k-c_i(\pi)$, where $h_k$ denotes the height of the $k$-th
step $p_k$ in the ballot path $p_1p_2\cdots p_k$.

\begin{figure}[h,t]
\begin{center}
\begin{picture}(380,140)
\setlength{\unitlength}{1mm}

\put(20,5){$\Longrightarrow $}

\put(33,0){\circle*{1.5}}\put(33,0){\line(1,1){5}}\put(34,1.8){\small0}
\put(38,5){\circle*{1.5}}\put(38,5){\line(1,1){5}}\put(39,6.8){\small1}
\put(43,10){\circle*{1.5}}\put(43,10){\line(1,1){5}}\put(44,11.8){\small1}
\put(48,15){\circle*{1.5}}\put(48,15){\line(1,-1){5}}\put(50.5,12){\small0}
\put(53,10){\circle*{1.5}}\put(53,10){\line(1,-1){5}}\put(55.7,6.8){\small1}
\put(58,5){\circle*{1.5}}\put(58,5){\line(1,1){5}}\put(59,7){\small1}
\put(63,10){\circle*{1.5}}

\put(70,5){$\Longrightarrow $}

\put(82,0){\circle*{1.5}}\put(82,0){\line(1,1){5}}\put(83,1.8){\small0}
\put(87,5){\circle*{1.5}}\put(87,5){\line(1,1){5}}\put(88,6.8){\small1}
\put(92,10){\circle*{1.5}}\put(92,10){\line(1,1){5}}\put(93,11.8){\small1}
\put(97,15){\circle*{1.5}}\put(97,15){\line(1,-1){5}}\put(99,11.8){\small0}
\put(102,10){\circle*{1.5}}\put(102,10){\line(1,-1){5}}\put(104,6.8){\small1}
\put(107,5){\circle*{1.5}}\put(107,5){\line(1,1){5}}\put(108.5,7.2){\small1}
\put(112,10){\circle*{1.5}}\put(112,10){\line(1,1){5}}\put(114,11.8){\small2}
\put(117,15){\circle*{1.5}}

\put(0,30){\circle*{1.5}}\put(0,30){\line(1,1){5}}\put(1,31.8){\small0}
\put(5,35){\circle*{1.5}}

\put(10,32){$\Longrightarrow $}
\put(20,30){\circle*{1.5}}\put(20,30){\line(1,1){5}}\put(21,31.8){\small0}
\put(25,35){\circle*{1.5}}\put(25,35){\line(1,1){5}}\put(26,36.8){\small1}
\put(30,40){\circle*{1.5}}

\put(35,32){$\Longrightarrow $}
\put(45,30){\circle*{1.5}}\put(45,30){\line(1,1){5}}\put(46,31.8){\small0}
\put(50,35){\circle*{1.5}}\put(50,35){\line(1,1){5}}\put(51,36.8){\small1}
\put(55,40){\circle*{1.5}}\put(55,40){\line(1,1){5}}\put(56,41.8){\small1}
\put(60,45){\circle*{1.5}}

\put(65,32){$\Longrightarrow $}
\put(75,30){\circle*{1.5}}\put(75,30){\line(1,1){5}}\put(76,31.8){\small0}
\put(80,35){\circle*{1.5}}\put(80,35){\line(1,1){5}}\put(81,36.8){\small1}
\put(85,40){\circle*{1.5}}\put(85,40){\line(1,1){5}}\put(86,41.8){\small1}
\put(90,45){\circle*{1.5}}\put(90,45){\line(1,-1){5}}\put(92,42){\small0}
\put(95,40){\circle*{1.5}}

\put(100,32){$\Longrightarrow $}
\put(110,30){\circle*{1.5}}\put(110,30){\line(1,1){5}}\put(111,31.8){\small0}
\put(115,35){\circle*{1.5}}\put(115,35){\line(1,1){5}}\put(116,36.8){\small1}
\put(120,40){\circle*{1.5}}\put(120,40){\line(1,1){5}}\put(121,41.8){\small1}
\put(125,45){\circle*{1.5}}\put(125,45){\line(1,-1){5}}\put(127,41.8){\small0}
\put(130,40){\circle*{1.5}}\put(130,40){\line(1,-1){5}}\put(132.5,36.9){\small1}
\put(135,35){\circle*{1.5}}
\end{picture}

\caption{The construction of $\Phi(\pi)$ for $\pi=2\bar
1547\bar 6\bar 3$.}\label{f1}
\end{center}
\end{figure}
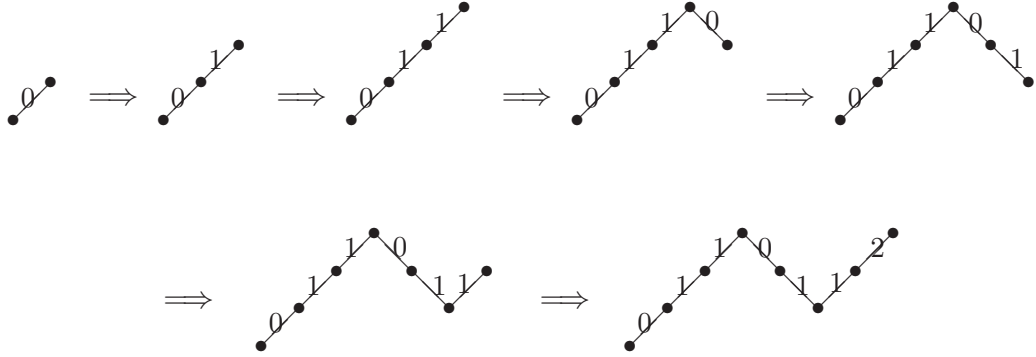

For example, let $n=7$ and
 $\pi=2\bar 1547\bar 6\bar 3$.
The construction of $\Phi(\pi)$ is illustrated  in Figure
\ref{f1}.

We now turn to the
 inverse map $\Psi$ from a labeled  ballot path
$(P;W)=(p_1p_2\cdots p_n;$  
$w_1w_2\cdots w_n)$ to a snake
$\pi=\pi_1\pi_2\cdots \pi_n$ of type $B_n$.

We shall construct a sequence of permutations
$\Gamma_0,\Gamma_{1},\Gamma_2,\ldots,\Gamma_n$, such that
$\Gamma_0=\emptyset$ and $\Gamma_n=\pi$ is the desired
 snake  of
type $B_n$. To reach this goal, we generate
 a sequence of labeled ballot
paths $ (P_1;W_1),(P_2;W_2),\cdots,(P_{n-1};W_{n-1})$,
 where $(P_1;W_1)=(P;W)$, and $P_{i+1}$ is
obtained from $P_{i}$ by contracting
 a certain step $p_{r_i}$ of
$P_{i}$ into a single point, and
 $W_{i+1}$ is obtained from $W_i$
by deleting the label of the step $p_{r_i}$ and updating the
labels of the other steps. Notice that
 $P_i$ has $n-i+1$ steps
 and $W_i$ has $n-i+1$ elements for $1\leq i\leq n$.
Below is a procedure to determine
 $(P_{i+1};W_{i+1})$ and $\Gamma_i$ from
  $(P_{i};W_i)$ and $\Gamma_{i-1}$.
Let us consider two cases.

\noindent Case 1: $P_i$ has an  odd  number of steps.

  If there  exists
   a down step in $P_i$ whose
  label equals its height, then we assume that
   $p_{r_i}$ is the leftmost among such
    down steps. Contract  $p_{r_i}$
  into a single point to form a ballot path $P_{i+1}$ and
  add 1 to  the labels of all down steps of $P_{i+1}$. Let
   $(P_{i+1};W_{i+1})$ denote the resulting labeled ballot path and
   set $\Gamma_i=\overline {n-r_i+1}
   \Gamma_{i-1}$.

  For the case that the label of any down step
   $P_i$  is less than its height,
    as will be shown, there must exist
 at least one
  up step labeled by 0. We assume that
  $p_{r_i}$ is the rightmost among such up steps. Contract
   $p_{r_i}$
  into a single point to form a ballot path $P_{i+1}$. Then subtract 1
  from the labels of up steps of $P_{i+1}$
 that are originally to the right of $p_{r_i}$  and
 add 1 to the labels of down steps of $P_{i+1}$
 that are originally to the left of $p_{r_i}$.
  Denote the resulting labeled ballot
 path by $(P_{i+1};W_{i+1})$ and set $\Gamma_i=
(n-r_i+1)\Gamma_{i-1}$.

\noindent Case 2:   $P_i$ has an even number of steps.

  If there
exists
   a down step of $P_i$ whose
  label equals 0, we assume that $p_{r_i}$ is
   the leftmost
   among such  down steps. Contract $p_{r_i}$
  into a single point to form a ballot path $P_{i+1}$.
   Then add 1
  to the labels of up steps of $P_{i+1}$
   which are originally to
   the right of $p_{r_i}$ and subtract  1 from
   the labels of down steps of $P_{i+1}$
   which are originally to the left of $p_{r_i}$.
    Denote the
   resulting labeled ballot path by $(P_{i+1};W_{i+1})$ and set
    $\Gamma_i=(n-r_i+1)\Gamma_{i-1}$.

For the case that there are no down steps in $P_i$ labeled by 0,
 as can be seen, there must exist at least one up step whose label equals its
height. We assume that
 $p_{r_i}$ is
the rightmost among such up steps. Contract $p_{r_i}$
  into a single point to form a ballot path $P_{i+1}$. Then subtract 1
  from the labels of all down steps of $P_{i+1}$.
Denote the resulting  path by $(P_{i+1};W_{i+1})$ and set
$\Gamma_i=\overline {n-r_i+1}\Gamma_{i-1}$.

\begin{figure}[h,t]
\begin{center}
\begin{picture}(380,150)
\setlength{\unitlength}{1mm}

\put(115,35){$\Longrightarrow $}

\put(82,30){\circle*{1.5}}\put(82,30){\line(1,1){5}}\put(83,31.8){\small0}
\put(87,35){\circle*{1.5}}\put(87,35){\line(1,1){5}}\put(88,36.8){\small1}
\put(92,40){\circle*{1.5}}\put(92,40){\line(1,1){5}}\put(93,41.8){\small1}
\put(97,45){\circle*{1.5}}\put(97,45){\line(1,-1){5}}\put(99,42){\small1}
\put(102,40){\circle*{1.5}}\put(102,40){\line(1,1){5}}\put(103.2,42){\small1}
\put(107,45){\circle*{1.5}}\put(107,45){\line(1,1){5}}\put(108,47){\small2}
\put(112,50){\circle*{1.5}} \put(82,22){$P_2=p_1p_2p_3p_4p_6p_7$}

\put(65,35){$\Longrightarrow $}

\put(23,30){\circle*{1.5}}\put(23,30){\line(1,1){5}}\put(24,31.8){\small0}
\put(28,35){\circle*{1.5}}\put(28,35){\line(1,1){5}}\put(29,36.8){\small1}
\put(33,40){\circle*{1.5}}\put(33,40){\line(1,1){5}}\put(34,41.8){\small1}
\put(38,45){\circle*{1.5}}\put(38,45){\line(1,-1){5}}\put(40,42){\small0}
\put(43,40){\circle*{1.5}}\put(43,40){\line(1,-1){5}}\put(45,37){\small1}
\put(48,35){\circle*{1.5}}\put(48,35){\line(1,1){5}}\put(49,36.8){\small1}
\put(53,40){\circle*{1.5}}\put(53,40){\line(1,1){5}}\put(54,41.7){\small2}
\put(58,45){\circle*{1.5}} \put(23,22){$P_1=p_1p_2p_3p_4p_5p_6p_7$}


\put(28,2){$\Longrightarrow $}
\put(0,0){\circle*{1.5}}\put(0,0){\line(1,1){5}}\put(1,1.8){\small0}
\put(5,5){\circle*{1.5}}\put(5,5){\line(1,1){5}}\put(6,6.8){\small1}
\put(10,10){\circle*{1.5}}\put(10,10){\line(1,-1){5}}\put(12,7){\small0}
\put(15,5){\circle*{1.5}}\put(15,5){\line(1,1){5}}\put(16,6.8){\small1}
\put(20,10){\circle*{1.5}}\put(20,10){\line(1,1){5}}\put(21.5,11.9){\small2}
\put(25,15){\circle*{1.5}} \put(-2,-8){$P_3=p_1p_3p_4p_6p_7$}

\put(63,2){$\Longrightarrow $}
\put(40,0){\circle*{1.5}}\put(40,0){\line(1,1){5}}\put(41,1.8){\small0}
\put(45,5){\circle*{1.5}}\put(45,5){\line(1,-1){5}}\put(47,2){\small0}
\put(50,0){\circle*{1.5}}\put(50,0){\line(1,1){5}}\put(51,1.8){\small0}
\put(55,5){\circle*{1.5}}\put(55,5){\line(1,1){5}}\put(57,7){\small1}
\put(60,10){\circle*{1.5}} \put(37,-8){$P_4=p_3p_4p_6p_7$}

\put(94,2){$\Longrightarrow $}
\put(75,0){\circle*{1.5}}\put(75,0){\line(1,1){5}}\put(76,1.8){\small0}
\put(80,5){\circle*{1.5}}\put(80,5){\line(1,1){5}}\put(81,6.8){\small1}
\put(85,10){\circle*{1.5}}\put(85,10){\line(1,1){5}}\put(86,11.8){\small2}
\put(90,15){\circle*{1.5}} \put(70,-8){$P_5=p_3p_6p_7$}

\put(120,2){$\Longrightarrow $}
\put(105,0){\circle*{1.5}}\put(105,0){\line(1,1){5}}\put(106,1.8){\small0}
\put(110,5){\circle*{1.5}}\put(110,5){\line(1,1){5}}\put(111,6.8){\small1}
\put(115,10){\circle*{1.5}} \put(100,-8){$P_6=p_6p_7$}


\put(130,0){\circle*{1.5}}\put(130,0){\line(1,1){5}}\put(131,1.8){\small0}
\put(135,5){\circle*{1.5}} \put(125,-8){$P_7=p_6$}
\end{picture}
\vspace*{1cm} \caption{The construction of  $\Psi(P;W)$ for the
labeled ballot path in Figure \ref{ff1}.}\label{fg}
\end{center}
\end{figure}
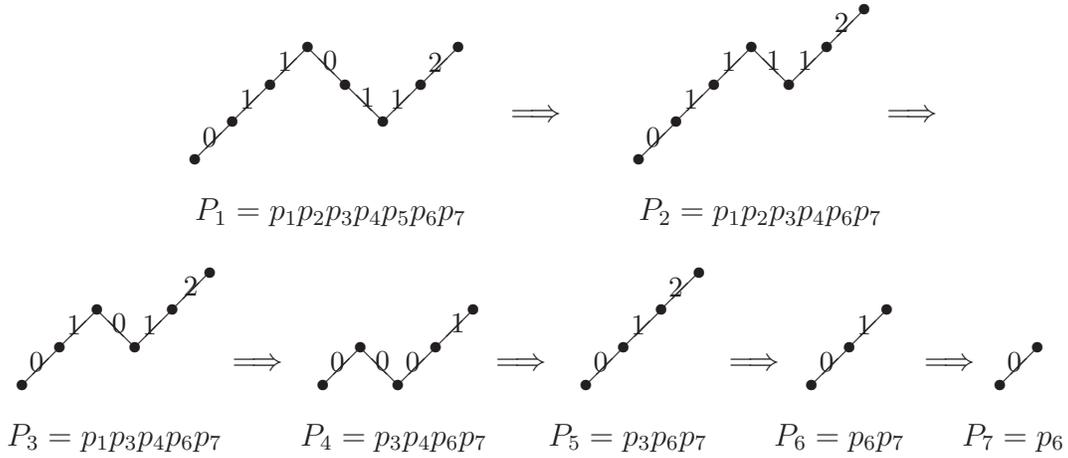

For the labeled ballot path $(P;W)=(uuudduu;0110112)$
  in Figure \ref{ff1},
  the construction of  $\Psi(P;W)$ is shown in Figure
\ref{fg}. The indices of the  steps $p_{r_i}$
 that are contracted are listed below:
$r_1=5,r_2=2,r_3=1,r_4=4,r_5=3,r_6=7,r_7=6$. The labeled ballot
paths $(P_i;W_i)$ are given in Figure \ref{fg}, and the
permutations $\Gamma_i$ are given as follows:
\[ \Gamma_0=\emptyset,\Gamma_1=\bar 3,
\Gamma_2=\bar6\bar3,\Gamma_3=7\bar6\bar3,\Gamma_4=47\bar6\bar3,
\Gamma_5=547\bar6\bar3,\Gamma_6=\bar 1547\bar6\bar3,
\Gamma_7=2\bar1547\bar6\bar3.\]

\section{The proof}

In this section, we shall show that the map $\Phi$ described in
the previous section is indeed a bijection.

\begin{theo}\label{t33}
The map $\Phi$ is a bijection between labeled ballot paths of
length $n$ and snakes of type $B_n$.
\end{theo}

\noindent{\it Proof. }
 As the first step, we verify
that $\Phi$ is well-defined, that is,
 for any snake $\pi=\pi_1\pi_2\cdots\pi_n$ of
type $B_n$, $\Phi(\pi)$ is a  labeled ballot paths.

Before we show that
 $\Phi(\pi)=(P;W)$
 is a labeled ballot path, it is necessary to
 prove that $P=p_1p_2\cdots p_n$ is a ballot path, that is,
  for
any $1\leq k\leq n$, the number of up steps is not less
 than the number
of down steps among the first $k$ steps of $P$.
 By the definition of
$\Phi$, we have $p_1=u$. Assume that in Step $k $ in the
implementation of $\Phi$, we have already constructed
$p_1,p_2,\ldots,p_{k-1}$ which form a ballot path.
 The
task of this step is to locate $n-k+1$ or $\overline {n-k+1}$ in
$\pi$ in order to get $p_k$. We consider two cases.

If $\pi_i=n-k+1$ and $i$ is odd or $\pi_i=\overline{n-k+1}$ and
$i$ is even, then we set $p_k=u$. Clearly, $p_1p_2\cdots p_k$ is a
ballot path.  Otherwise, we set $p_k=d$ and we wish to show that
the height $h_k$ of $p_k$ is nonnegative.
 Consider the case
 $\pi_i=\overline {n-k+1}$
 and $i$ is odd. Observe that the
height of $p_k$
 is the number of up steps among
 $p_1,p_2,\cdots,p_{k-1}$ subtracts the number of down
steps  among
 $p_1,p_2,\cdots,p_{k}$. By
the definition of $\Phi$, we have
\begin{align}
h_k&=\#\{1\leq j\leq n\,|\,\pi_j>0, n-k+1<\pi_j\ {\rm and} \ j \
{\rm is\
odd}\}\nonumber\\
&\qquad  +\#\{1\leq j\leq n\,|\,\pi_j<0, n-k+1<|\pi_j|\ {\rm and} \ j \ {\rm is\ even}\}\nonumber\\
&\qquad  -\#\{1\leq j\leq n\,|\, \pi_j<0,n-k+1\leq |\pi_j|\ {\rm and} \ j \ {\rm is\ odd}\}\nonumber\\
&\qquad   - \#\{1\leq j\leq n\,|\,\pi_j>0, n-k+1<\pi_j\ {\rm and}\
j\ {\rm is\ even}\}.\label{down2}
\end{align}
 In view of the alternating property of
$\pi$, if there is a negative element $\pi_{2i+1}$ at an odd
position of $\pi$, then $\pi_{2i}$   must be negative as well and
$\pi_{2i+1}>\pi_{2i}$. Consequently,
\begin{align*}
& \#\{1\leq j\leq n\,|\, \pi_j<0,n-k+1\leq |\pi_j|\ {\rm and} \ j \ {\rm is\ odd}\} \\
&\quad \leq \#\{1\leq j\leq n\,|\,\pi_j<0, n-k+1<|\pi_j|\ {\rm and}
\ j \ {\rm is\ even}\}.
\end{align*}
On the other hand, if there is
 a positive element $\pi_{2i}$
 at an even position of $\pi$, then $\pi_{2i-1}$
 must be positive as well and $\pi_{2i}<\pi_{2i-1}$. This yields
 that
\begin{align*}
& \#\{1\leq j\leq n\,|\,\pi_j>0, n-k+1<\pi_j\ {\rm and}\
j\ {\rm is\ even}\} \\
&\quad \leq \#\{1\leq j\leq n\,|\,\pi_j>0, n-k+1<\pi_j\ {\rm and} \
j \ {\rm is\ odd}\}.
\end{align*}
 Thus we deduce that
 whenever there is a negative term contributing to
$h_k$, there is at least one positive term. So we conclude that
$h_k\geq 0$. A similar argument applies to
 the case that $\pi_i=n-k+1$ and $i$ is even.  Hence we have shown that $P$ is a ballot
path.

 We next prove that the label of any step in $\Phi(\pi)$
  is nonnegative and it does not
  exceed its
height. Let $\pi=\pi_1\cdots\pi_i\cdots\pi_n$. Assume we are in
the Step $k$ and we have determined $(p_1\ldots
p_{k-1};w_1,\ldots,w_{k-1})$, which is a labeled ballot path of
length $k-1$. We proceed to
 locate $n-k+1$ or $\overline {n-k+1}$
in $\pi$ in order to determine $p_k$ and its label $w_k$. Suppose
that  $\pi_i=\overline {n-k+1}$
 and  $i$ is   odd.
 In this case, by the definition of $\Phi$,   we have $p_k=d$
  and  $w_k=h_k-c_i(\pi)$.
 We claim
that  $c_i(\pi)\leq h_k$.  In computing $h_k$ by using formula
(\ref{down2}), we shall split the range of $j$ into two cases: one
case is $1\leq j\leq i$ and the other case is $i+1\leq j\leq n$.
In other words, we shall consider the contributions of
$\pi_1\pi_2\ldots \pi_{i}$ and $\pi_{i+1}\ldots \pi_n$ to the
value of $h_k$.

We claim that $c_i(\pi)$ is less than or equal to the contribution
of $\pi_{i+1}\ldots\pi_{n}$ to $h_k$.
 Suppose that $n$ is odd. By
the definition of $c_i(\pi)$, a pair $(\pi_{2j},\pi_{2j+1})$ of
consecutive elements of $\pi$ with
 $i<2j\leq n-1$  contributes 1
to the value of $c_i(\pi)$  if $\pi_{2j}<\pi_i<\pi_{2j+1}<0$ or
$\pi_{2j}<\pi_i<0$ and $\pi_{2j+1}>0$. If there is a pair
$(\pi_{2j},\pi_{2j+1})$  with $\pi_{2j}<\pi_i<\pi_{2j+1}<0$, then
this pair contributes 1 to both $h_k$ and $c_i(\pi)$. If there is a
pair $(\pi_{2j},\pi_{2j+1})$  with $\pi_{2j}<\pi_i<0$ and
$\pi_{2j+1}>0$, then this pair contributes  1 or 2 to $h_k$ (depends
on whether $|\pi_{2j+1}|$ is greater than $n-k+1$), while
contributes exactly 1 to $c_i(\pi)$. It is straightforward to check
that if a pair $(\pi_{2j},\pi_{2j+1})$ does not contribute to
$c_i(\pi)$, then it contributes 0 or 1 to $h_k$. On the hand hand,
because $\pi_1\ldots\pi_{i}$ contributes 0 to $c_i(\pi)$, it remains
to show that the contribution of $\pi_1\ldots\pi_{i}$ to $h_k$ is
nonnegative. Let
\begin{align*}
g_i(\pi)&=\#\{1\leq j\leq i\,|\,\pi_j>0, n-k+1<\pi_j\ {\rm and} \ j \ {\rm is\ odd}\}\\
&\qquad  +\#\{1\leq j\leq i\,|\,\pi_j<0, n-k+1<|\pi_j|\ {\rm and} \ j \ {\rm is\ even}\}\\
& \qquad  -\#\{1\leq j\leq i\,|\, \pi_j<0,n-k+1\leq |\pi_j|\ {\rm and} \ j \ {\rm is\ odd}\}\\
&\qquad   - \#\{1\leq j\leq i\,|\,\pi_j>0, n-k+1<\pi_j\ {\rm and}\
j\ {\rm is\ even}\}.
\end{align*}
 By the same reasoning as in the
 proof of  $h_k\geq 0$,
  we can verify  that $g_i(\pi)\geq0$.
 Thus we have completed the proof for the case that
$\pi_i=\overline {n-k+1}$ and both $n,i$ are odd. All the other
cases depending on the sign of $\pi_i$ and the parities of $n$ and
$i$ can be treated in the same manner. Hence the details are
omitted.

Our next task is to show that the map $\Psi$ is well-defined,
namely, for any labeled ballot path $(P;W)$ of length $n$, the
signed permutation
 $\pi=\Gamma_{n}=\pi_1\pi_2\cdots\pi_n$ is a
snake of
 type
$B_n$, i.e., $0<\pi_1>\pi_2<\pi_3>\cdots \pi_n$.

Suppose that at the $i$-th step we have already constructed a
labeled ballot path $(P_{i};W_i)$. We first consider the case that
 $P_i$ has an odd number of steps. In this case
 we aim to  show that after
contracting a certain step of $P_i$,
 we can get a ballot path
$P_{i+1}$. By our construction of $\Psi$, if there is
a down step
in $P_i$ whose label
 equals its height, then we
contract the leftmost such down step in $P_i$. In this case,
  we automatically get a ballot path $P_{i+1}$.
  Otherwise, we consider the case that
   there exist no such down steps.
    In particular, this implies that there are
 no down steps with height 0.
 By our construction $\Psi$, we shall
 contract the rightmost up step labeled by 0.
  After we contract
 this up step, it is easily seen that
 every step in $P_{i+1}$ has nonnegative height
 since we know that there are no down steps that touch
 the $x$-axis. So we also get a ballot path $P_{i+1}$
  in this case.
One can check that after we update the labels of the steps in
$P_{i+1}$, each step  will  have a nonnegative label that is less
than or equal to its height.

The case that $P_i$ has an even number of steps can be dealt with
by the same argument as for the case that $P_i$ has an odd number
of steps. Thus we conclude that once we have accomplished the
mission in step $i$, we are led to a labeled ballot path
$(P_{i+1}, W_{i+1})$ and a signed permutation
 $\Gamma_{i}$.

 We now turn to the proof of the alternating property
 of $\pi$. It is apparent from the construction
 of $\Psi$ that $\pi_1>0$. Now we prove that
 $\pi_1>\pi_2<\pi_3>\ldots \pi_n$. Suppose that in step $i-1$ we have already
  constructed a
 signed permutation $\Gamma_{i-1}$ and a labeled ballot
 path $(P_{i};W_{i})$.
To determine $\Gamma_i$, by our construction,
 we are supposed to contract a certain step $p_{r_{i}}$
 in $P_{i}$
 to form a ballot path $P_{i+1}$ and to
 set $\Gamma_i=(n-r_i+1)\Gamma_{i-1}$
 or $\Gamma_i=\overline{n-r_i+1}\Gamma_{i-1}$ depending on
 whether $p_{r_{i}}$ is an up step or a down step.

To determine $\Gamma_{i+1}$, by our construction,
 we are supposed to contract a certain step
  $p_{r_{i+1}}$ of $P_{i+1}$
 to form a ballot path $P_{i+2}$ and to set
  $\Gamma_{i+1}=(n-r_{i+1}+1)
 \Gamma_{i}$
 or $\Gamma_{i+1}=\overline{n-r_{i+1}+1}
 \Gamma_{i}$ depending on whether
 $p_{r_{i+1}}$ is an up step or a down step.
For notational convenience, set $t_i=n-r_i+1$ and $
t_{i+1}=n-r_{i+1}+1$. There are four possibilities
 for the
construction of $\Gamma_{i+1}$, namely,
 $t_{i+1}t_i\Gamma_{i-1},
 \bar t_{i+1}\bar t_i\Gamma_{i-1}$,
$\bar t_{i+1} t_i\Gamma_{i-1}$ and  $t_{i+1}\bar t_i\Gamma_{i-1}$.

We only consider the
 case that $P_i$ has an odd number
of steps and so $t_i$ is at an odd position of $\pi$. To prove the
alternating
 property of $\Gamma_n$, it is necessary to verify
 that $t_{i+1}<t_{i}$,  $\bar
t_{i+1}<\bar t_i$ and $\bar t_{i+1}<t_i$. And the situation that
$\Gamma_{i+1} =t_{i+1}\bar t_i\Gamma_{i-1}$ can never happen.

 In this case, in the $i$-th step,
  suppose that we contract a down step $p_{r_i}$ of $P_i$,
 and in
the $(i+1)$-st step,  suppose that we contract an up step
$p_{r_{i+1}}$ of $P_{i+1}$. By the construction of $\Psi$, we have
$\Gamma_{i+1}=\bar t_{i+1}\bar t_i\Gamma_{i-1}$. We claim that
$t_{i}<t_{i+1}$, i.e., $r_i>r_{i+1}$. Otherwise, we may assume that
$r_i<r_{i+1}$. Once the down step $p_{r_i}$ is contracted, the
height of all steps to the right of the step $p_{r_i}$
 will
increase by 1, but by the construction of $\Psi$,
the labels of up steps remain unchanged. This implies
that  the labels of up steps to the
right of $p_{r_i}$ cannot be equal to
 their heights. Therefore, the up step
$p_{r_{i+1}}$ cannot be chosen in the $(i+1)$-st step,
which is
 a contradiction. So we deduce that
 $\bar t_{i+1}<\bar t_i$. The discussions for the cases that $t_{i+1}<t_{i}$,
  $\bar t_{i+1}<t_i$ and the situation that
$\Gamma_{i+1} =t_{i+1}\bar t_i\Gamma_{i-1}$ can never happen are
similar.

We are now left with the case that the number of steps of $P_i$ is
even to complete the proof of the alternating property. But the
argument in this case is analogous to that for the case that $n-i$
is even. Hence we have reached the conclusion that
 $\Gamma_{n}=\pi_1\pi_2\cdots\pi_n$ is a
snake of
 type
$B_n$.

Finally, we wish to confirm that the maps $\Phi$ and $\Psi$ are
 inverses of each other.
 Because both  $\Phi$ and $\Psi$ are carried out
  in $n$ steps, it suffices to verify that
    the $i$-th step of $\Phi$ and the $i$-th step of $\Psi$
    are inverses of each other.
 Let
 $\pi=\pi_1\pi_2\cdots\pi_n$. For $1\leq i\leq n$, let
 $\Pi_i=\pi_1\pi_2\cdots\pi_i$. Define $\Phi(\Pi_i)$ to be the
 labeled ballot path by applying $\Phi$
 to the standardization of
 $\Pi_i$. The standardization of a snake $\Pi_i$ is a snake obtained by keeping the
 sign of each element
 unchanged and replacing the smallest
 element   by 1, the second smallest element
 by 2, and so forth.
 Note that one can
 apply $\Phi$ to $\Pi_1, \cdots,\Pi_n$ step by step and
 finally obtain $\Phi(\pi)=\Phi(\Pi_n)$.
 By the construction of $\Psi$, one can check that
 $\Phi(\Pi_i)=(P_{n-i+1};W_{n-i+1})$ for $1\leq i\leq n$.
 That is, the inverse procedure to derive $\Phi(\Pi_{i+1})$
 from $\Phi(\Pi_i)$
 coincides with the procedure to construct $\Psi(P_{n-i+2};W_{n-i+2})$
 from $\Psi(P_{n-i+1};W_{n-i+1})$.
 Therefore $\Phi$ and
 $\Psi$ are inverses of each other.

 In summary, we have shown that that
 the map $\Phi$ is a bijection between labeled ballot paths of
length $n$ and snakes of type $B_n$. \qed

\section {A refinement}

In this section, we obtain the bivariante generating function for
the number $B(n,k)$ of labeled ballot paths of length $n$ that end
at a given point $(n,k)$, where $0\leq k\leq n$. The numbers
$B(n,k)$ can be considered as a refinement of the Springer numbers.
By restriction of the bijection $\Psi$, we also obtain a bijection
between labeled Dyck paths of length $2n$ and  alternating
permutations on $[2n]$. By considering the last step of a labeled
ballot path, it is easy to derive  the following recurrence
relation.

\begin{theo} For $1\leq k \leq n$, we have
\begin{align}\label{Bnk}
B(n,k)=(k+1)B(n-1,k+1)+kB(n-1,k-1).
\end{align}
\end{theo}

Note that in the above recurrence relation we need the convention
that $B(n,k)=0$ for $n<k$. Moreover, since a ballot path can never
ends at a point $(m,n)$ where $m+n$ is odd, so $B(n,k)=0$ if $n+k$
is odd.
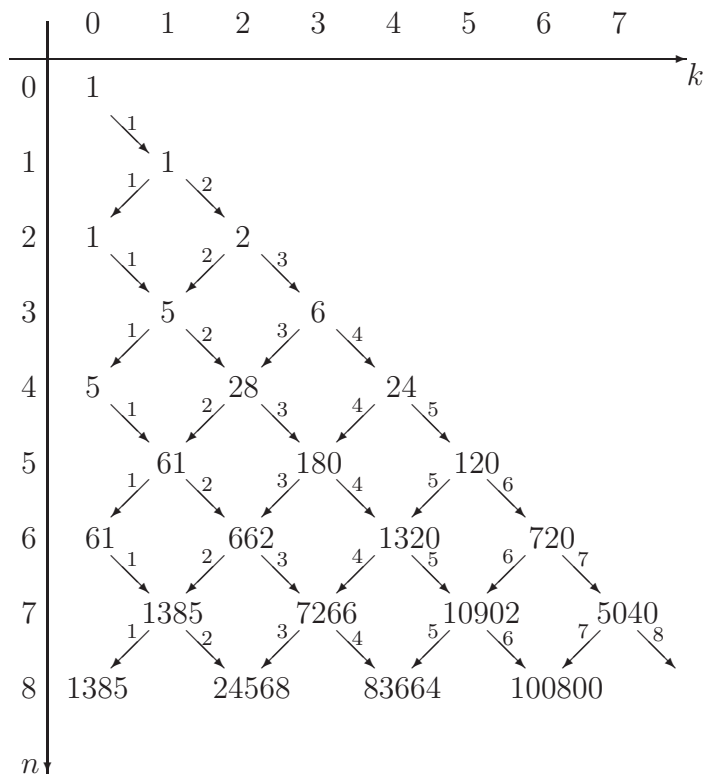
\begin{figure}[h]
\hspace*{4cm}
\begin{picture}(15,40)
\setlength{\unitlength}{0.5cm}

\put(0,1){\vector(0,-1){20}}\put(-1,0){\vector(1,0){18}}
\put(17,-0.7){$k$}\put(-0.7,-19){$n$}
\put(1,0.7){0}\put(3,0.7){1}\put(5,0.7){2}\put(7,0.7){3}\put(9,0.7){4}\put(11,0.7){5}\put(13,0.7){6}\put(15,0.7){7}
\put(-0.7,-1){0}\put(-0.7,-3){1}\put(-0.7,-5){2}\put(-0.7,-7){3}\put(-0.7,-9){4}\put(-0.7,-11){5}\put(-0.7,-13){6}\put(-0.7,-15){7}\put(-0.7,-17){8}
\put(1,-1){1}
\put(1,-5){1}\put(1,-9){5}\put(1,-13){61}\put(0.5,-17){1385}
\put(1.7,-1.5){\vector(1,-1){1}}\put(2.1,-1.9){\scriptsize 1}
\put(2.7,-3.2){\vector(-1,-1){1}}\put(2.1,-3.4){\scriptsize 1}
\put(1.7,-5.2){\vector(1,-1){1}}\put(2.1,-5.5){\scriptsize 1}
\put(2.7,-7.2){\vector(-1,-1){1}}\put(2.1,-7.4){\scriptsize 1}
\put(1.7,-9.2){\vector(1,-1){1}}\put(2.1,-9.5){\scriptsize 1}
\put(2.7,-11.2){\vector(-1,-1){1}}\put(2.1,-11.4){\scriptsize 1}
\put(1.7,-13.2){\vector(1,-1){1}}\put(2.1,-13.5){\scriptsize 1}
\put(2.7,-15.2){\vector(-1,-1){1}}\put(2.1,-15.4){\scriptsize 1}

\put(3,-3){1}\put(3,-7){5}\put(2.9,-11){61}\put(2.5,-15){1385}
\put(3.7,-3.2){\vector(1,-1){1}}\put(4.1,-3.5){\scriptsize 2}
\put(4.7,-5.2){\vector(-1,-1){1}}\put(4.1,-5.4){\scriptsize 2}
\put(3.7,-7.2){\vector(1,-1){1}}\put(4.1,-7.5){\scriptsize 2}
\put(4.7,-9.2){\vector(-1,-1){1}}\put(4.1,-9.4){\scriptsize 2}
\put(3.7,-11.2){\vector(1,-1){1}}\put(4.1,-11.5){\scriptsize 2}
\put(4.7,-13.2){\vector(-1,-1){1}}\put(4.1,-13.4){\scriptsize 2}
\put(3.7,-15.2){\vector(1,-1){1}}\put(4.1,-15.6){\scriptsize 2}

\put(5,-5){2}\put(4.8,-9){28}\put(4.8,-13){662}\put(4.4,-17){24568}
\put(5.7,-5.2){\vector(1,-1){1}}\put(6.1,-5.5){\scriptsize 3}
\put(6.7,-7.2){\vector(-1,-1){1}}\put(6.1,-7.4){\scriptsize 3}
\put(5.7,-9.2){\vector(1,-1){1}}\put(6.1,-9.5){\scriptsize 3}
\put(6.7,-11.2){\vector(-1,-1){1}}\put(6.1,-11.4){\scriptsize 3}
\put(5.7,-13.2){\vector(1,-1){1}}\put(6.1,-13.5){\scriptsize 3}
\put(6.7,-15.2){\vector(-1,-1){1}}\put(6.1,-15.4){\scriptsize 3}

\put(7,-7){6}\put(6.6,-11){180}\put(6.6,-15){7266}
\put(7.7,-7.2){\vector(1,-1){1}}\put(8.1,-7.5){\scriptsize 4}
\put(8.7,-9.2){\vector(-1,-1){1}}\put(8.1,-9.4){\scriptsize 4}
\put(7.7,-11.2){\vector(1,-1){1}}\put(8.1,-11.5){\scriptsize 4}
\put(8.7,-13.2){\vector(-1,-1){1}}\put(8.1,-13.4){\scriptsize 4}
\put(7.7,-15.2){\vector(1,-1){1}}\put(8.1,-15.6){\scriptsize 4}

\put(9,-9){24}\put(8.8,-13){1320}\put(8.4,-17){83664}
\put(9.7,-9.2){\vector(1,-1){1}}\put(10.1,-9.5){\scriptsize 5}
\put(10.7,-11.2){\vector(-1,-1){1}}\put(10.1,-11.4){\scriptsize 5}
\put(9.7,-13.2){\vector(1,-1){1}}\put(10.1,-13.5){\scriptsize 5}
\put(10.7,-15.2){\vector(-1,-1){1}}\put(10.1,-15.4){\scriptsize 5}

\put(10.8,-11){120}\put(10.5,-15){10902}
\put(11.7,-11.2){\vector(1,-1){1}}\put(12.1,-11.5){\scriptsize 6}
\put(12.7,-13.2){\vector(-1,-1){1}}\put(12.1,-13.4){\scriptsize 6}
\put(11.7,-15.2){\vector(1,-1){1}}\put(12.1,-15.6){\scriptsize 6}

\put(12.8,-13){720}\put(12.3,-17){100800}
\put(13.7,-13.2){\vector(1,-1){1}}\put(14.1,-13.5){\scriptsize 7}
\put(14.7,-15.2){\vector(-1,-1){1}}\put(14.1,-15.4){\scriptsize 7}

\put(14.6,-15){5040}
\put(15.7,-15.2){\vector(1,-1){1}}\put(16.1,-15.5){\scriptsize 8}
\end{picture}
\vspace*{10cm} \caption{The recurrence relation for $B(n,k)$}\label{TB}
\end{figure}

Note that when $k=0$,  $B(2n,0)$ is the number of
labeled Dyck paths
of length $2n$, where a labeled Dyck path of length $2n$
is a labeled ballot path of length $2n$ that ends with
a point on the $x$-axis. It is worth mentioning that the numbers
$B(2n,0)$  are in fact   the secant numbers and
they are closely related to alternate
level codes of ballots, see Strehl \cite{Strehl}.
Recall that an
alternate level code of ballots of length $n$ is an
 integer sequence
$\lambda=\lambda_1\lambda_2\cdots\lambda_n$ such that $\lambda_1=1$, and for $2\leq j\leq
n$,
\[ \lambda_{j-1}+1\geq\lambda_j\geq 1 .\]
Denote by $\Lambda_n$ the set of  alternate level codes of ballots of length
$n$.  For example,
\[ \Lambda_3=\{111,112,121,122,123\}.\]

Rosen \cite{Rosen} derived the following formula
\begin{align}\label{tan}
\sum_{n\geq0}\left(\sum_{\lambda\in\Lambda_n}\prod_{i=1}^n\lambda_i(\lambda_i+1)\right)
{x^n\over n!}=\tan x.
\end{align}
 Strehl \cite{Strehl} deduced
 the  secant companion equation of (\ref{tan}):
\begin{align}\label{sec}
\sum_{n\geq0}\left(\sum_{\lambda\in\Lambda_n}\prod_{i=1}^n\lambda^2_i\right)
{x^n\over n!}=\sec x.
\end{align}

To make a connection between labeled
Dyck paths and alternate level codes of ballots, we need the following bijection, see Stanley \cite[Ex. 6.19]{Stanley}.

\begin{theo}
There is a bijection between the set of Dyck paths of length $2n$ and the
set of  alternate level codes of ballots of length $n$.
\end{theo}

\proof Let $\lambda=\lambda_1\lambda_2\cdots\lambda_n\in\Lambda_n$ be an alternate level
 code of ballots  of length $n$.
  For convenience, we set $\lambda_{n+1}=1$.
  We shall  construct
 a Dyck path $P$ of length $2n$ from
 $\lambda$.
 Let $P=P_1P_2\cdots P_n$, where $P_i=ud^k$ for some
 $k\geq0$, that is, $P_i$ consists of an up step followed by $k$ down steps.
 For $1\leq i\leq n$,  if $\lambda_i=\lambda_{i+1}-1$ then $k=0$, $P_i=u$.
 If $\lambda_i\geq \lambda_{i+1}$ then
  $P_i=ud^{\lambda_i-\lambda_{i+1}+1}$.
 It is necessary to show
 that the above construction generates
   a Dyck path of length $2n$. That is,
   after the $i$-th step,
    the number of down steps is less than or equal to
    the number of up steps. That is, we wish to show that
 \[\sum_{j=1}^i(\lambda_j-\lambda_{j+1}+1)\leq i.\]
Since $\lambda_1=1$ and $\lambda_{i+1}\geq 1$,
the above inequality
is apparently true. Moreover,
 it is  easy to check that there are $n$ down
steps, namely,
\[\sum_{j=1}^n(\lambda_j-\lambda_{j+1}+1)=n.\]
Conversely, given a Dyck path of length $2n$,
 let $\lambda_i$ be
the height of the
 $i$-th up step plus one.
 It is readily verified
  that  $\lambda=\lambda_1\lambda_2\cdots\lambda_n$
 is  an alternate level
 code of ballots of length $n$. This completes the proof.
\qed

For instance, let $\lambda=122\in\Lambda_3$. Then the Dyck path
corresponding to $\lambda$ is $uududd$. Using the above
bijection£¬ we are led to a connection between the number
$B(2n,0)$ and  alternate level codes of ballots.

\begin{coro} We have
\begin{align}\label{alt}
B(2n,0)=\sum_{\lambda\in\Lambda_n}\prod_{i=1}^n\lambda^2_i.
\end{align}
\end{coro}

\proof Relation (\ref{alt}) follows from the
observation that for a given Dyck
path, the number of labelings equals the product of
 squares of the elements of the corresponding  alternate level code of ballots. \qed

 In passing, we mention that
 Getu, Shapiro and Woen \cite{Getu}
 have considered a generalization
  of the formula  of Rosen \cite{Rosen}
  on tangent numbers, namely,   equation (\ref{tan}).
  More precisely, for a
 given ballot path, they defined
 the weight of the path to be the
 product of the $y$-coordinate  of all the endpoints,
  except for the
 last point.  Let
 $T(n,k)$ denote the sum of   weights of  ballot paths from
 $(1,1)$ to $(n,k)$. It is easily checked that
\[T(n,k)=(k-1)T(n-1,k-1)+(k+1)T(n-1,k+1).\]
   When
 $k=1$, $T(n,1)$ is the tangent number, that is,
 \[\sum_{n\geq1}T(n,1){x^n\over n!}=\tan x. \]
 They gave a table for $T(n,k)$ similar to the table
  in Figure \ref{TB},
 where the first column
 consists of the tangent numbers. For $k\geq 1$,
 they obtained the generating function
\[\sum_{n\geq 1}T(n,k){x^n\over n!}={\tan^k x\over k}.\]

 By replacing the first column of their table
  by the secant numbers they introduced another number
   $E(n,k)$, and
  they considered the following  recurrence relation
\[E(n,k)=(k -1)E(n-1,
k -1)+ kE(n-1, k+1),\]
 where $E(0,1)=1, E(1,2)=E(2,1)=1$ and $E(n,k)=0$
 for $n<k-1$ or
 $k<1$. When $k=1$, $E(n,1)$ is the secant number.
 However, no  combinatorial interpretation was given
 for the numbers $E(n,k)$.
Using the recurrence relation of $E(n,k)$,
 Getu, Shapiro and Woen
\cite{Getu}  derived the exponential generating function
\begin{equation}\label{Fnk}
\sum\limits_{n\geq
k}E(n,k)\frac{x^n}{n!}=\tan^{k-1} x\sec x.
\end{equation}
Comparing the recurrence relations and
  initial values of $B(n,k)$ and $E(n,k)$, it became
  apparent that
 \[B(n,k)=E(n,k+1).\]
 Therefore $B(n,k)$ can be viewed as a
 combinatorial explanation for $E(n,k)$.
 Moreover we obtain the generating functions $G_n(y)$ for
the rows of the table for $B(n,k)$.
Let
\begin{equation}\label{Gn}
G_n(y)=\sum\limits_{0\leq k\leq n}B(n,k)y^k.
\end{equation}
Note that
\[G_n(1)=\sum\limits_{0\leq k\leq n}B(n,k)\]
equals the $n$-th Springer number.
 Let $B(x,y)$ be the generating
function for $G_n(y)$, that is,
\[
B(x,y)=\sum\limits_{n\geq0}G_n(y){x^n\over n!}.
\]
Then we have the following formula.

\begin{theo}
\begin{equation} \label{bxy}
B(x,y)=\frac{1}{\cos x-y\sin x}.
\end{equation}
\end{theo}

\proof Let \[F_k(x)=\sum_{n\geq k}B(n,k){x^n\over n!}=\sum_{n\geq
k}E(n,k+1){x^n\over n!}=\tan^k x\sec x.\] Therefore,
\begin{align*}
B(x,y)&=\sum_{ n\geq 0}\sum_{0\leq k\leq n}B(n,k)y^k{x^n\over
n!}=\sum_{ k\geq 0}F_k(x)y^k=\frac{1}{\cos x -y\sin x},
\end{align*}
as required.\qed

To conclude this paper, we give two  applications of the bijection
$\Phi$. More precisely, we obtain
 a classification of snakes of
type $B_n$, and we establish a
connection between  labeled Dyck paths and
alternating permutations.

Define the following statistic
\begin{align*} \alpha(\pi)&=\#\{1\leq j\leq n\,|\,\pi_j>0\ {\rm and} \ j
\ {\rm is\ odd}\}\\
&\qquad+\#\{1\leq j\leq n\,|\,\pi_j<0\ {\rm and} \ j \ {\rm is\ even}\}\\
&\qquad  -\#\{1\leq j\leq n\,|\, \pi_j<0\ {\rm and} \ j \ {\rm is\
odd}\}\\
&\qquad - \#\{1\leq j\leq n\,|\,\pi_j>0 \ {\rm and}\ j\ {\rm is\
even}\}.
\end{align*}
Then we have  the following
 classification of snakes of type $B_n$.
\begin{theo}\label{t4}
For $0\leq k \leq n$, $B(n,k)$ equals
the number of snakes $\pi=\pi_1\pi_2\cdots\pi_n$ with
$\alpha(\pi)=k$.
\end{theo}

In particular, let us consider the
implication of the above theorem for $k=0$.
Recall that $B(2n,0)$ is the number of labeled Dyck paths of
length $2n$. By (\ref{sec}) and (\ref{alt}), we have
\[\sum_{n\geq 0}B(2n,0){x^n\over n!}=\sec x.\]
 It now comes to our mind that $\sec x$ is the
 generating function for the number $E_{2n}$
 of alternating
 permutations on $[2n]$. This indicates that
 $B(2n,0)$ equals $E_{2n}$. The following theorem asserts
 that the restriction of the bijection $\Phi$ to
 labeled Dyck paths serves as a combinatorial
 interpretation of the fact that $B(2n,0)=E_{2n}$.
 Roughly speaking, when restricted to
 labeled Dyck paths the map $\Psi$ does
  not involve
 any negative elements and when restricted to alternating
 permutations the map $\Phi$ generates labeled Dyck paths.

 The following theorem is concerned with the restriction
 of the map $\Psi$. It is not difficult to see that the
 restriction of $\Psi$ to labeled Dyck path is the inverse
 of the restriction of $\Phi$ to alternating permutations.

\begin{theo}
The map $\Psi$ induces a bijection between
labeled Dyck paths of
length $2n$ and alternating permutations on $[2n]$.
\end{theo}

\proof Let $(P;W)=(p_1\ldots p_{2n};w_1\ldots w_{2n})$
 be a labeled
Dyck path of length $2n$. We wish to show that
$\pi=\Psi(P;W)=\pi_1\ldots\pi_{2n}$
contains no negative elements. In the
first step of  $\Psi$, since $(P;W)$ is a labeled
Dyck path, there must exist down steps labeled by 0.
 Assume that
$p_{r_1}$ is the leftmost among such down steps.
Applying the map $\Psi$, we are supposed to contract
$p_{r_1}$ into a single point to form a ballot path $P_2$.
Then we are supposed to add
1 to the labels of up steps of $P_2$ which are originally
 to the
right of $p_{r_1}$ and subtract 1 from the labels
 of down steps of
$P_2$ which are originally to the left of $p_{r_1}$.
Hence we get a
labeled ballot path $(P_2;W_2)$ and a permutation
$\Gamma_1=(n-r_1+1)\Gamma_0=(n-r_1+1)$, which contains no
negative elements.

Similarly, in Step 2, in the labeled ballot path $(P_2;W_2)$, there
does not exist any down step of $P_2$ whose label equals its height.
So we can find an up step of $P_2$ labeled by 0. Suppose that
$p_{r_2}$ is the leftmost up step  of $P_2$ with label 0.
Contracting $p_{r_2}$ gives  a ballot path $P_3$. Then subtract 1
from the labels of up steps of $P_3$ that are originally to the
right of $p_{r_2}$ and add 1 to the labels of down steps of $P_3$
that are originally to the left of $p_{r_2}$. Thus we obtain a
labeled ballot path $(P_3;W_3)$ and a permutation
$\Gamma_2=(n-r_2+1)\Gamma_1=(n-r_2+1)(n-r_1+1)$ without negative
elements.

Since $P_1$ is a Dyck path of length $2n$, and an up step in $P_1$
and a down step in $P_2$ are contracted,
   there are $n-1$ up steps and $n-1$ down steps in $P_3$.
   It follows that
 $(P_3;W_3)$ is a labeled Dyck path.
 Continuing the above process, we eventually
 obtain an alternating permutation.

 Conversely, given   an alternating permutation
  $\pi=\pi_1\pi_2\ldots
 \pi_{2n}$  of length $2n$, we wish to show
 that $\Phi(\pi)=p_1\ldots p_{2n}$ is a labeled Dyck path
 of length $2n$. Since
 $\Phi(\pi)$ is a labeled ballot path already, it is enough to show
 that it has the same number of up steps as down steps.
  In  Step $k$ of the map $\Phi$,
  we are supposed to find the location of
  the element $n-k+1$ in
  $\pi$. Assume that $\pi_i=n-k+1$.
   Carrying out the construction of $\Phi(\pi)$,
   this is what happens: if $i$ is odd, then we have
  $p_k=u$, and if $i$ is even, then we have $p_k=d$.
  Since $\pi$ has $2n$ elements, so we conclude with
   $n$ up steps as well as $n$ down
  steps.
 Therefore $\Phi(\pi)$ is a labeled
 Dyck path. This completes the proof. \qed

\vspace{0.5cm}
 \noindent{\bf Acknowledgments.}  This work was supported by  the 973
Project, the PCSIRT Project of the Ministry of Education,  and the National Science Foundation of China.

\end{document}